\documentclass[11pt]{article}
\usepackage{amssymb}
\usepackage{amsmath}

\newcommand{\dar}{{\downarrow}}

\newtheorem{proposition}{Proposition}[section]

\newtheorem{corollary}[proposition]{Corollary}
\newtheorem{definition}[proposition]{Definition}
\newtheorem{theorem}[proposition]{Theorem}

\newcounter{opgaveteller}
\newcounter{lijst-teller}

\newcounter{boon}
\newcounter{boon2}
\newenvironment{proof}{\noindent {\bf Proof}. \nopagebreak }{\nopagebreak\hfill\rule{2mm}{3mm}}

\input xypic
\xyoption{all}
\UseComputerModernTips
\CompileMatrices
\definemorphism{dat}\dashed\tip\notip
\definemorphism{dub}\Solid\Tip\notip

\newenvironment{rlist}%
   {\begin{list}{\roman{lijst-teller}\/)\hfil}%
              {\labelwidth 2em%
               \leftmargin\labelwidth\advance\leftmargin by\labelsep%
               \usecounter{lijst-teller}}}%
   {\end{list}}

   {\end{list}}

    {\end{list}}

\title{A general form of relative recursion}
\author{Jaap van Oosten\\ Department of Mathematics\\ Utrecht University\\ P.O.Box 
80.010, 3508 TA Utrecht, The Netherlands\\ {\tt jvoosten@math.uu.nl}}
\begin{document}
\maketitle

\begin{abstract} 
\noindent The purpose of this note is to observe a generalization of the concept 
``computable in\ldots '' to arbitrary partial combinatory algebras. For every partial 
combinatory algebra (pca) $A$ and every partial endofunction on $A$, a pca $A[f]$ is 
constructed such that in $A[f]$, the function $f$ is representable by an element; a 
universal property of the construction is formulated in terms of Longley's 2-category of 
pcas and decidable applicative morphisms.

It is proved that there is always a geometric inclusion from the realizability topos on 
$A[f]$ into the one on $A$, and that there is a meaningful preorder on the partial 
endofunctions on $A$ which generalizes Turing reducibility.\end{abstract}


\section*{Introduction}
In \cite{LongleyJ:reatls}, John Longley defined a 2-category of partial combinatory algebras
(see~\ref{pcasec} and \ref{2catsec} for definitions). The morphisms are different from what one might expect:
rather that `algebraic' maps, they are more like simulations (of one world of computation in another).
Accordingly, a morphism from $A$ to $B$ is a total relation between the underlying sets.

Longley's definition made a lot of sense since there are nice functorial connections between pcas
and their corresponding realizability categories (realizability toposes and categories of assemblies).

However, the 2-category has not been studied in great detail. It does not appear to have a lot
of categorical structure, and not much is known. Fundamental questions, such as: which properties
of partial combinatory algebras are stable under isomorphism, or equivalence?, have not been answered
(indeed, such questions have hardly been posed).

In this paper, I present a simple construction which is available in this category: {\em adjoin a partial function}.
That is, given a pca $A$ and a partial endofunction $f$ on $A$, construct a pca $A[f]$ in which the function $f$
is `computable'. $A[f]$ should, of course, possess a universal property, and this property is formulated
with respect to what Longley calls `decidable' morphisms.

Characteristically for the non-algebraic flavour of the 2-category, $A[f]$ is not constructed by adding elements,
but by modifying the application function. We obtain results generalizing the situation of computing relative to
an oracle: a preorder, similar to (and generalizing) Turing reducibility, can be defined on the partial endofunctions
on $A$; and there is always a geometric inclusion from the realizability topos on $A[f]$ into the one on $A$.

It is also a surprising corollary of this work that every total pca is isomorphic to a nontotal one.
\subsection{Basic notions and notations}
\subsubsection{Partial combinatory algebras}\label{pcasec}
A {\em partial combinatory algebra\/} (pca) is a set $A$ together with a partial function 
$A\times A\rightharpoonup A$ called {\em application}, which satisfies a few conditions. 
We write the application as $(a,b)\mapsto ab$ or $a{\cdot}b$. $ab\dar$ means that the 
application $ab$ is defined. When dealing with compound terms like $(ac)(bc)$, the 
definedness of the term is meant to imply the definedness of every subterm. For terms $t$ 
and $s$, the notation $t\simeq s$ means that $t$ is defined exactly when $s$ is; and that 
they deonte the same element when defined. $t=s$ will mean $t\simeq s$ and $t\dar$. As 
usual, we associate to the left, that is: $abc$ means $(ab)c$. Elements of $A$ are 
usually called {\em combinators}.

With these conventions, $(A,{\cdot})$ is a pca iff there are combinators $K$ and $S$ in 
$A$ satisfying, for all $a,b,c\in A$:\begin{itemize}
\item $Kab=a$
\item $Sab\dar$
\item $Sabc\simeq ac(bc)$\end{itemize}
For a careful account of the theory of pcas, see \cite{BethkeI:npca} or 
\cite{LongleyJ:reatls}. We recall a few properties.

In a pca $A$ there is a choice of Booleans $\top$ and $\bot$, and a `definition by cases' 
combinator $C$ such that for all $a,b\in A$, $C\top ab=a$ and $C\bot ab=b$; $C$ is 
pronounced (and written) as {\sf If}\ldots {\sf then}\ldots {\sf else}\ldots .

In $A$ there is a choice of elements $\overline{n}$ for every natural number $n$, such 
that for every partial recursive function $F$ of $k$ variables there is a combinator 
$a_F$ such that for every $k$-tuple $(n_1,\ldots ,n_k)$, 
$a_F\overline{n_1}\cdots\overline{n_k}\dar$ precisely when $F(n_1,\ldots ,n_k)$ is 
defined, and $a_F\overline{n_1}\cdots\overline{n_k}=\overline{F(n_1,\ldots ,n_k)}$ if 
this is the case. There is a coding of finite sequences of elements of $A$, together with 
combinators which allow us to manipulate them: if we write $[u_0,\ldots ,u_{n-1}]$ for 
the code of the sequence $(u_0,\ldots ,u_{n-1})$, there is a combinator {\sf lh} which 
gives the {\em length\/} of the coded sequence (i.e. ${\sf lh}[u_0,\ldots 
,u_{n-1}]=\overline{n}$), there are combinators picking the $i$-th element of the coded 
sequence (we simply write $u_i$ for its effect) and a concatenation operator; we write 
$[u_0,\ldots u_{n-1}]\ast [v_0,\ldots 
,v_{m-1}]$ for the effect of this last combinator.

All these facts follows from the existence, in $A$, of a combinator for primitive 
recursion. Moreover, in every pca $A$ there is a fixpoint combinator $Y$ satisfying: 
$Yf\dar$ for all $f\in A$, and $Yfa\simeq f(Yf)a$. We shall refer to this fact as `the 
recursion theorem in $A$'.

Every pca $A$ is `combinatory complete': for every term $t$ (constructed from variables, 
constants from $A$, and the application function) and every sequence of variables 
$x_1,\ldots x_{n+1}$ which contains all variables in $t$, there is an element $\Lambda 
^*x_1\cdots x_{n+1}.t$ in $A$ which satisfies for all $a_1,\ldots a_{n+1}$ in 
$A$:\begin{itemize}
\item $(\Lambda ^*x_1\cdots x_{n+1}.t)a_1\cdots a_n\dar$
\item $(\Lambda ^*x_1\cdots x_{n+1}.t)a_1\cdots a_{n+1}\simeq t(a_1,\ldots 
,a_{n+1})$\end{itemize}

\subsubsection{Longley's 2-category of pcas; assemblies; decidable maps}\label{2catsec}
The following definition is due to John Longley (\cite{LongleyJ:reatls}).
\begin{definition}\label{apmor}\em Let $A$ and $B$ be pcas. An {\em applicative 
morphism\/} from $A$ to $B$ is a function $\gamma$ from $A$ to the set ${\cal P}^*(B)$ of 
nonempty subsets of $B$, such that there exists an element $r\in B$ with the property 
that if $aa'\dar$ in $A$, $b\in\gamma (a)$ and $b'\in\gamma (a')$, then $rbb'\dar$ and 
$rbb'\in\gamma (aa')$. The element $r$ is said to be a {\em realizer\/} for $\gamma$.\end{definition}
Given two applicative morphisms $\gamma :A\to B$ and $\delta :B\to C$, the composition 
$\delta\gamma :A\to C$ is the function $a\mapsto\bigcup_{b\in\gamma (a)}\delta (b)$ from 
$A$ to ${\cal P}^*(C)$. It is easy, using combinatory completeness, to find a realizer for
$\delta\gamma$ in terms of realizers for $\gamma$ and $\delta$.

This composition is evidently associative and has identities 
$a\mapsto\{ a\}$, so we have a category of pcas.

This category is preorder-enriched: given two applicative morphisms $\gamma ,\delta :A\to 
B$, we say $\gamma\preceq\delta$ if there is an $s\in B$ such that for all $a\in A$ and 
all $b\in\gamma (a)$, $sb\in\delta (a)$. We say that $\gamma$ and $\delta$ are isomorphic 
if $\gamma\preceq\delta$ and $\delta\preceq\gamma$ both hold.

Two pcas are {\em equivalent\/} if there are $\gamma :A\to B$ and $\delta :B\to A$ such 
that both composites are isomorphic to identities.
\medskip

\noindent An {\em assembly\/} on a pca $A$ is a set $X$ together with a map $E_X: X\to 
{\cal P}^*(A)$. If $(X,E_X)$ and $(Y,E_Y)$ are assemblies on $A$, a map of assemblies is 
a function $f:X\to Y$ such that there is an element $r\in A$ such that for all $x\in X$ 
and all $a\in E_X(x)$, $ra\dar$ and $ra\in E_Y(f(x))$. One says that the element $r$ {\em 
tracks\/} the function $f$. Assemblies on $A$ and maps of assemblies form a category 
${\rm Asm}(A)$. This category is regular and comes equipped with an adjunction to the 
category of Sets: the forgetful (or global sections) functor $\Gamma :{\rm Asm}(A)\to 
{\rm Set}$ is left adjoint to the functor $\nabla :{\rm Set}\to {\rm Asm}(A)$ which sends 
a set $X$ to the pair $(X,E_X)$ where $E_X(x)=A$ for all $x\in X$.

An important justification for definition~\ref{apmor} is the following theorem by 
Longley: every applicative morphism $\gamma :A\to B$ determines a regular functor $\gamma 
^*:{\rm Asm}(A)\to {\rm Asm}(B)$ which commutes with the functors $\Gamma$; conversely, 
every such functor is induced by an applicative morphism which is unique up to 
isomorphism.

Note, that $\gamma :A\to B$ establishes $A$ as an assembly on $B$.
\begin{definition}\label{decmor}\em A morphism $\gamma :A\to B$ is {\em decidable\/} if 
there is an element $d\in B$ (the {\em decider\/} for $\gamma$) such that if $\top 
_A,\bot _A$ are the Booleans in $A$ and $\top _B,\bot _B$ the Booleans in $B$, for every 
$b\in\gamma (\top _A)$ we have $db=\top _B$ and for every $b\in\gamma (\bot _A)$, 
$db=\bot _B$.\end{definition}
In \cite{LongleyJ:reatls}, Longley proved
\begin{proposition}\label{deccop}An applicative morphism $\gamma :A\to B$ is decidable if 
and only if the corresponding functor $\gamma ^*:{\rm Asm}(A)\to {\rm Asm}(B)$ preserves 
finite coproducts. Moreover this is equivalent to: $\gamma ^*$ preserves the natural 
numbers object.\end{proposition}
\begin{corollary}\label{deccor}If $\delta =\gamma\zeta$ is a commutative triangle of applicative morphisms
such that $\delta$ and $\zeta$ are decidable, then so is $\gamma$.\end{corollary}

\section{Definition of $A[f]$ and basic properties}
\begin{definition}\em Let $\gamma :A\to B$ be an applicative morphism of pcas and $f:A\rightharpoonup A$ 
a partial function. We say that $f$ is {\em representable\/} w.r.t.\ $\gamma$ if there is an
element $r_f\in B$ such that for every $a\in {\rm dom}(f)$ and every $b\in\gamma (a)$, $r_fb\dar$
and $r_fb\in\gamma (f(a))$. We say that $f$ is representable in $A$ if $f$ is representable w.r.t.\
the identity morphism on $A$.\end{definition}
The representability of $f$ with respect to $\gamma$ can also be seen as follows: let $({\rm dom}(f),\gamma )$ 
be the regular sub-assembly of $(A,\gamma )$ (as assemblies on $B$). Then $f$ is representable with
respect to $\gamma$ if and only if $f$ is a map of assemblies: $({\rm dom}(f),\gamma )\to (A,\gamma )$.
\begin{theorem}\label{main}For every pca $A$ and every partial endofunction $f$ on $A$ there exist
a pca $A[f]$ and a decidable applicative morphism $\iota _f:A\to A[f]$ with the following properties:
\begin{rlist}\item $f$ is representable w.r.t.\ $\iota _f$;
\item for every decidable applicative morphism $\gamma :A\to B$ such that $f$ is representable w.r.t.\
$\gamma$, there is a decidable applicative morphism $\gamma _f:A[f]\to B$ such that $\gamma _f\iota _f=\gamma$,
and $\gamma _f$ is unique with this property. Moreover, if $\delta: A[f]\to B$ is such that $\delta\iota _f\cong\gamma$, 
then $\delta\cong\gamma _f$\end{rlist}\end{theorem}
\begin{proof} For the construction of $A[f]$, let's agree on some notation for codes of finite sequences:
if $u=[u_o,\ldots ,u_{n-1}]$ and $i<n$, $u^{<i}$ denotes $[u_0,\ldots ,u_{i-1}]$ and $u^{\geq i}$ denotes
$[u_i,\ldots ,u_{n-1}]$; for $i\leq j<n$, $^{i\leq}u^{<j}$
denotes $[u_i,\ldots ,u_{j-1}]$. Let $p,p_0,p_1$ be pairing and projection combinators in $A$, i.e. satisfying
for all $a,b\in A$: $p_0(pab)=a$ and $p_1(pab)=b$. Let {\sf Not} be a combinator such that ${\sf Not}\top =\bot$ and
${\sf Not}\bot =\top$. 

The underlying set of $A[f]$ will be $A$. We define a new application ${\cdot}^f$ on $A$ as follows.
For $a,b\in A$, an $f$-{\em dialogue\/} between $a$ and $b$ is a code of a sequence $u=[u_0,\ldots ,u_{n-1}]$ 
such that for all $i<n$ there is a $v_i\in A$ such that
$$a{\cdot}([b]\ast u^{<i})=p\bot v_i\mbox{ and }f(v_i)=u_i$$
We say that $a{\cdot}^fb$ is defined with value $c$, if there is an $f$-dialogue $u$ between $a$ and $b$ 
such that $$a{\cdot}([b]\ast u)=p\top c$$
We show first, that $(A,{\cdot}^f)$ is a pca.

Let $K_f=\Lambda ^*x.p\top (\Lambda ^*y.p\top x_0)$. Then clearly $K_f{\cdot}^fa=\Lambda ^*y.p\top a$ for all 
$a\in A$, so $(K_f{\cdot}^fa){\cdot}^fb=a$ for all $a,b\in A$.

For the combinator $S_f$, by primitive recursion it is possible to construct a term $t(x,y)$ of $A$ such
that for all $u$, the application $t(x,y){\cdot}u$ is given by the following instructions:

\noindent $t(x,y){\cdot}u=$\begin{quote}$xu$ if $\forall i\leq {\sf lh}u\,\,{\sf Not}(p_0(xu^{<i}))$.

If $i$ is minimal such that $p_0(xu^{<i})$, let $\alpha =p_1(xu^{<i})$ and output $y([u_0]\ast u^{\geq i})$ if
$\forall j(i\leq j<{\sf lh}u\to {\sf Not}p_0(y([u_0]\ast ^{i\leq}u^{<j}))$.

If $j$ is minimal such that $p_0(y([u_0]\ast ^{i\leq}u^{<j}))$, let $\beta = p_1(y([u_0]\ast ^{i\leq}u^{<j}))$
and output $\alpha ([\beta ]\ast u^{\geq j})$ if
$\forall k(j\leq k<{\sf lh}u\to {\sf Not}(p_0(\alpha ([\beta ]\ast ^{j\leq}u^{<k}))))$.

If $k$ is minimal such that $(p_0(\alpha ([\beta ]\ast ^{j\leq}u^{<k})))$, output $(p_1(\alpha ([\beta ]\ast
^{j\leq}u^{<k})))$.\end{quote}
Note, that $t(a,b){\cdot}^fc\simeq (a{\cdot}^fc){\cdot}^f(b{\cdot}^fc)$ for all $a,b,c$. Therefore,
let $$S_f=\Lambda ^*x.p\top (\Lambda ^*y.p\top t(x_0,y_0))$$ 
Then $(S_f{\cdot}^fa){\cdot}^fb=t(a,b)$ for all $a$
and $b$. This establishes $A[f]$ as a pca.

Note that the combinators $K_f$ and $S_f$ don't really depend on $f$. This is analogous to the fact that for
a coding of Turing machine computations with oracle $U$, the $S^m_n$-functions are primitive recursive, and do not
depend on $U$.

The map $\iota _f:A\to A[f]$ given by $a\mapsto\{ a\}$ is an applicative morphism $A\to A[f]$. Indeed, if $ab=c$
then $(\Lambda ^*x.p\top (ax_0)){\cdot}^fb=c$; so if $r=\Lambda ^*yx.p\top (y_0x_0)$ then $r$ realizes $\iota _f$.

The decidability of $\iota _f$ is left to the reader.

For the universal property, suppose $\gamma :A\to B$ is a decidable applicative morphism which is realized by $r$ 
and let $d$ be a decider for $\gamma$. Moreover suppose that $\overline{f}$ represents $f$ w.r.t.\ $\gamma$.

Let $\pi _0,\pi _1\in B$ be such that if $b\in\gamma (a)$ then $\pi _ib\in\gamma (p_ia)$. Similarly, let $C$ and $C'$
in $B$ be such that if $b\in\gamma (a)$ and $v\in\gamma (u)$ then $Cbv\in\gamma ([a]\ast u)$ and 
$C'bv\in\gamma (u\ast [a])$.

Now use the recursion theorem in $B$ to find an element $U$ such that for all $b,b',v$:
$$\begin{array}{rcl}Ubb'v & \simeq & \mbox{{\sf If} } d(\pi _0(rb(Cb'v))) \\
 & & \mbox{{\sf then} }\pi _1(rb(Cb'v)) \\
 & & \mbox{{\sf else} }Ubb'(C'(\overline{f}(\pi _1(rb(Cb'v))))v) \end{array}$$
The reader can check the following: suppose $u$ is an $f$-dialogue between $a$ and $a'$ in $A$, 
$b\in\gamma (a),b'\in\gamma (a')$, $i<{\sf lh}u$, $v\in\gamma (u^{<i})$ and $w=C'(\overline{f}(\pi _1(rb(Cb'v))))v$.
Then $w\in\gamma (u^{\leq i})$ and $Ubb'v=Ubb'w$. Furthermore, if $u$ is such that $a([a']\ast u)=p\top c$, then 
$Ubb'v\in\gamma
(c)$.

Therefore, choose $e\in\gamma ([\, ])$ and let
$$\rho =\Lambda ^*xx'.Uxx'e$$
Then $\rho$ realizes $\gamma$ as applicative morphism: $A[f]\to B$. We denote this last morphism by $\gamma _f$.

Obviously, the diagram
$$\diagram A\rto^i\drto_{\gamma} & A[f]\dto^{\gamma _f} \\ & B\enddiagram$$
commutes on the nose. Moreover, since $\iota _f(a)=\{ a\}$, if $\delta :A[f]\to B$ were such that 
$\delta\iota\cong\gamma _f\iota$, then $\delta\cong\gamma _f$. So $\gamma _f$ is unique with respect to the property 
that the
diagram commutes on the nose, and essentially unique with respect to the property that it commutes up to isomorphism.
The decidability of $\gamma _f$ is a direct consequence of Corollary~\ref{deccor} and can also be verified directly. 
\end{proof}
\begin{corollary}\label{isocor}~~~~~~~~~~~~~~~~~~~~~~~~~~~~~~~~~~~~~~~~~~~~~~~\begin{rlist}
\item If $f$ is representable in $A$, then $A$ and $A[f]$ are isomorphic pcas.
\item If $f$ and $g$ are two partial endofunctions on $A$, the pcas $A[f][g]$ and $A[g][f]$ are isomorphic; we
may therefore write $A[f,g]$.
\item If ${\cal K}_1$ denotes Kleene's pca of partial recursive application, $f:\mathbb{N}\to\mathbb{N}$ is a
partial function and ${\cal K}_1^f$ is the pca of partial recursive application with an oracle for $f$, then
${\cal K}_1^f$ is isomorphic to ${\cal K}_1[f]$.
\item There exists a nontotal pca which is isomorphic to a total pca.\end{rlist}\end{corollary}
\begin{proof} The first two statements are immediate from the uniqueness statement in theorem~\ref{main}. The third
statement is easy. Finally, the fourth statement follows from the fact that $A[f]$ is never total (the element
$a=\Lambda ^*x.p\bot\bot$ is such that $a{\cdot}^fb$ is never defined), so if $A$ is total and $f$ is representable in $A$,
then $A\cong A[f]$ by i).\end{proof}
\medskip

\noindent {\bf Example} In \cite{OostenJ:casfft}, a total combinatory algebra $\cal B$ of partial functions on $\mathbb{N}$ is
defined, and it is proved that the representable functions are those functions which are continuous for the Scott
topology and satisfy some ``sequentiality'' condition. One might consider what happens if a ``parallel'' function is
adjoined to this: e.g.\ let $F:{\cal B}\to {\cal B}$ be the function such that for all $\alpha\in {\cal B}$, $F(\alpha )(0)=0$ 
if and only if $0\in {\rm dom}(\alpha )$ or $1\in {\rm dom}(\alpha )$ (and undefined else), and $F(\alpha )(n)$ is undefined for
all $n>0$. What would the representable functions of ${\cal B}[F]$ be? My conjecture would be that these are exactly
all Scott-continuous functions on $\cal B$.
\medskip

\noindent {\bf Remarks}\begin{enumerate}\item The construction of $A[f]$ induces a preorder on the set of partial endofunctions
of $A$, which generalizes Turing degrees: let $f\leq _A g$ if and only if $f$ is representable in $A[g]$ (with respect to $\iota
_g$). Since the diagram
$$\diagram A\rto\dto & A[g]\dto \\A[h]\rto & A[g,h]\enddiagram$$
commutes, it is easy to see that $\leq _A$ is a transitive relation (it is reflexive by ~\ref{main}(i)).
\item There is a universal solution to the problem of ``making $A$ decidable''; adjoin a function $f$ to $A$ where 
$$f(x)=\left\{\begin{array}{rl}\top & \mbox{if }p_0x=p_1x \\ \bot & \mbox{else}\end{array}\right.$$
\item This seems to be a good point to correct a claim made in \cite{HofstraP:ordpca}, lemma 5.4. It is claimed that no
total pca can be equivalent to a pca $A$ in which there is an element $z$ such that for all $x$, $zx\dar$ and $zx\neq x$.
However, this is established only if ``equivalent'' is replaced by ``isomorphic''. Therefore the original claim remains
an open problem. Another open problem, as far as I know, is: give an example of two pcas which are equivalent, but not
isomorphic.\end{enumerate}

\section{A geometric inclusion of realizability toposes}
The construction of $A[f]$ generalizes another aspect of relative recursion, known from the theory of realizability 
toposes. It is well known that for every pca $A$ there exists a {\em realizability topos\/} ${\rm RT}(A)$. The best
studied example is ${\rm RT}({\cal K}_1)$, the {\em effective topos}(\cite{HylandJ:efft}). In \cite{HylandJ:efft} and
\cite{PhoaW:rcet} it is explained that ${\rm RT}({\cal K}_1^f)$ is a {\em subtopos\/} of ${\rm RT}({\cal K}_1)$, in the
topos-theoretic sense. Here we shall see that this generalizes to geometric inclusions ${\rm RT}(A[f])\to {\rm RT}(A)$.

In \cite{HofstraP:ordpca}, the authors analyze a generalization of Longley's 2-category of pcas, and characterize which
applicative morphisms give rise to geometric morphisms between realizability toposes. The key concept is that of a {\em
computationally dense\/} morphism. Unfortunately, the definition given in l.c.\ is not quite adequate; see also
\cite{HofstraP:erropc}. I state the correct definition here for pcas.
\begin{definition}\label{compdense}\em Suppose that $F:A\to B$ is a function between pcas such that the map $a\mapsto\{ F(a)\}$ is
an applicative morphism. $F$ is {\em computationally dense\/} if there is an $m\in B$ with the property that for every
$b\in B$ one can find an $a\in A$ such that for all $a'\in A$:\begin{itemize}
\item[]If $bF(a')\dar$ in $B$, then $aa'\dar$ in $A$, and $mF(aa')=bF(a')$\end{itemize}\end{definition}
Let $P(A)$ and $P(B)$ denote the realizability triposes on $A$ and $B$. Then in \cite{HofstraP:ordpca} it is shown that
the map of indexed preorders induced by $F^*$ (where $F^*:{\cal P}(A)\to {\cal P}(B)$ sends $\alpha$ to $F[\alpha ]$) has
an indexed right adjoint if and only if $F$ is computationally dense.

In that case, the right adjoint is induced by the map $\hat{F}:{\cal P}(B)\to {\cal P}(A)$, given by
$$\hat{F}(\beta )=\{ a\in A\, |\, mF(a)\in\beta\}$$
where $m\in B$ witnesses the computational density of $F$.

It is easily verified then, that if $F$ is computationally dense and $m$ is as in definition~\ref{compdense}, then the geometric
morphism $(\hat{F},F^*)$ is an inclusion precisely when the following condition holds:
\begin{itemize}\item[(in)] There is a $c\in B$ such that for every $b\in B$ there is an $a\in A$ such that $cb=F(a)$ and $m(cb)=b$
\end{itemize}
\begin{proposition} The identity function $A\to A[f]$ is computationally dense and satisfies the condition {\rm
(in)}.\end{proposition}
\begin{proof} This is quite simple. Let $m$ be an element of $A$ such that for every $y\in A$ and every code of a sequence
$v$, $m([y]\ast v)\simeq yv$.

Given $b\in A$, let $a\in A$ be such that for all $a'\in A$, $aa'\simeq\Lambda ^*v.b([a']\ast v)$. Then $aa'$ is always
defined. Moreover,
$$m([aa']\ast v)\simeq (aa')v\simeq b([a']\ast v)$$
It follows that $m{\cdot}^f(aa')\simeq b{\cdot}^fa'$ in $A[f]$. This proves that the identity function is computationally dense.

Moreover, if $c=\Lambda ^*x.p\top (\Lambda ^*v.p\top v_0)$ then for all $a$, $c[a]=p\top (\Lambda v.p\top a)$; hence
$c{\cdot}^fa=\Lambda ^*v.p\top a$ and $$m([c{\cdot}^fa])=(c{\cdot}^fa)[\, ]=p\top a$$
so $m{\cdot}^f(c{\cdot}^fa)=a$, which proves (in).\end{proof}

\begin{small}
\bibliographystyle{plain}
 
\end{small}
\end{document}